\documentclass[12pt]{article}

\usepackage[english]{babel}
\usepackage{amsmath,amssymb,dsfont, natbib, amsthm,bbm}
\usepackage{setspace}
\usepackage{epsfig}
\usepackage{multirow}
\usepackage{booktabs}
\usepackage{xr}

\usepackage{color}

\graphicspath{{Figures/}}

\definecolor{purple}{rgb}{0.55,0.2,0.90}

\renewcommand{\hat}{\widehat}

\newcommand{\PP}{{\mathbb P}}

\newcommand{\RR}{\mathbb{R}}

\newcommand{\NN}{\mathbb{N}}

\newtheorem{theorem}{Theorem}[section]
\newtheorem{lemma}[theorem]{Lemma}

\title{Testing convexity of a discrete distribution}
\author{Fadoua Balabdaoui$^1$, C\'ecile Durot$^2$\thanks{Corresponding author. Email address: cecile.durot@gmail.com} and Fran\c{c}ois Koladjo$^{3}$}
\date{\vspace{-6ex}}

\begin{document}
\maketitle
\begin{center}
$^{1}$CEREMADE, Universit\'e Paris-Dauphine, 75775, Paris, France \\
$^{2}$Modal'x, Universit\'e Paris Nanterre, F-92001, Nanterre, France\\
$^{3}$ENSPD, Universit\'e Parakou, BP 55 Tchaourou,   B\'enin 
\end{center}

\begin{abstract}
Based on the convex least-squares estimator, we propose two different procedures for testing convexity of a probability mass function supported on $\NN$ with an unknown finite support. The procedures are shown to be asymptotically calibrated.
\end{abstract}

\section{The testing problem}\label{sec: pb}
Estimating a probability mass function (pmf) under a shape constraint has attracted attention in the very last years, see \cite{jankowski2009estimation}, \cite{durot2013least}, \cite{BJRP2013},  \cite{giguelay2016estimation}. With a more applied point of view, \cite{durot2015nonparametric} developped a nonparametric method for estimating the number of species under the assumption that the abundance distribution is convex. As the method applies only if the convexity assumption is fulfilled,  it would be sensible, before to implement it, to test whether or not the assumption is fulfilled. This motivates the present paper where a method for testing convexity of a pmf on $\NN$ is developped.

We consider i.i.d. observations $X_{1},\dots,X_{n}$  from an unknown pmf $p_{0}$  on $\NN$. Assuming that $p_{0}$ has a finite support $\{0,\dots,S\}$ with an unknown integer $S > 0$, 
we aim at testing the null hypothesis
$H_0$: "$p_{0}$ is convex on $\NN $" versus the alternative $H_1$: "$p_0$ is not convex."
With $p_n$ the empirical pmf (defined by 
$p_n(j) = {n}^{-1} \sum_{i=1}^n \mathbbm{1}_{\{X_i = j\}}$ for all $j \in \NN$), a natural procedure  rejects $H_0$ if $p_{n}$ is too far from ${\mathcal C}_{1}$, the set of all convex probability mass functions on $\NN$. 
Hence, with 
$\|q\|^2= \sum_{j \in \NN} (q(j))^2 $
for a sequence $q= \{q(j), j \in \NN\}$, 
we reject $H_0$ if $\inf_{p\in{\cal C}_{1}}\|p_{n}-p\|$ is too large. It is proved in \citet[Sections 2.1 to 2.3]{durot2013least} that the minimizer exists, is unique, and can be implemented with an appropriate algorithm, so our critical region takes the form
$\{ { T_n} >t_{\alpha,n}\}$
where $T_n=\sqrt n\| p_{n}-\hat p_{n}\|$, $\hat p_n$ is the minimizer of  
$ \|p_{n}-p\|^2$ 
over $p\in{\cal C}_{1}$, and $t_{\alpha,n}$ is an appropriate quantile. The main difficulty now is to find $t_{\alpha,n}$ in such a way that the corresponding test has asymptotic level $\alpha$.

We consider below two different constructions of $t_{\alpha,n}$. First, we will define $t_{\alpha,n}$ to be the $(1-\alpha)$-quantile of a random variable whose  limiting distribution coincides with the limiting distribution of  $T_n$ under $H_0$. Next, we will calibrate the test under the least favorable hypothesis. Both methods require knowledge of the limiting distribution of $T_n$ under $H_0$. To this end, we need notation. For all $p = \{p(j), j \in \NN\}$ and $k \in \NN\backslash\{0\}$ we set $\Delta p(k)=p(k+1)-2p(k)+p(k-1)$ (hence $p$ is convex on $\NN$ iff $\Delta p(k)\geq 0$ for all $k$) and a given $k\in\NN\backslash\{0\}$ is called  a knot of $p$ if $\Delta p(k)>0$. Furthermore, we 
denote by $g_0$  a $(S+2)$ centered Gaussian vector whose dispersion matrix $\Gamma_{0}$ has component $(i+1,j+1)$ equal to $\mathbbm{1}_{\{i= j\}} p_{0}(i)-p_{0}(i)p_{0}(j)$ for all  $i,j=0,\dots,S+1$, and by $\hat g_{0}$ the minimizer of 
$\sum_{k=0}^{S+1}  \left(g(k) - g_{0}(k)\right)^2$
over the set ${\cal K}_{0}$ of all functions $g=(g(0), \ldots, g(S+1)) \in\RR^{S+2} $ such that $\Delta g(k)\geq 0$ for all $k\in\{1,\dots,S\}$ with possible exceptions at the knots of $p_{0}$. Existence, uniqueness and characterization of $\hat g_{0}$ are given in \citet[Theorem 3.1]{BDK2014limit}. 
The asymptotic distribution of $T_n$ under $H_0$ is given below.

\begin{theorem}\label{cor: asympt}
Under $H_0$, $T_n\stackrel{d}{\longrightarrow}\hat T_0$ as $n\to\infty$, where $\hat T_0=\sum_{k=0}^{S+1}  \left(\hat g_0(k) - g_{0}(k)\right)^2$.
\end{theorem}

\section{Calibrating by estimating the limiting distribution}\label{sec: lim}
In order to approximate the distribution of $T_n$ under $H_0$, we will construct a random variable  that weakly converges to $\hat T_0$ (see Theorem \ref{cor: asympt}) and which can  be approximated via Monte-Carlo simulations.  To this end,  
let $S_{n}=\max\{X_{1},\dots,X_{n}\}$. Also, let $g_{n}$ be a random vector which, conditionally on $(X_{1},\dots,X_{n})$,  is distributed as a centered Gaussian vector of dimension $S_n+2$ with dispersion matrix $\Gamma_{n}$, the matrix with component $(i+1,j+1)$ equal to $\mathbbm{1}_{\{i=j\}}p_{n}(i)-p_{n}(i)p_{n}(j)$ for all $i,j=0,\dots,S_{n}+1$.  
Now, let $\hat g_{n}$ be the minimizer of 
$\sum_{k=0}^{S_{n}+1}  \left(g(k) - g_{n}(k)\right)^2$
over a set ${\cal K}_{n}$ that approaches ${\cal K}_{0}$ as $n\to \infty$. 
Below, we give an extended version  of \citet[Theorem 3.3]{BDK2014limit}, with the same choice for ${\cal K}_{n}$.

 \begin{theorem}\label{theo: CVcalib}
 Let $(v_{n})_{n\in\NN}$ be a sequence of positive numbers that satisfy
$ v_{n}=o(1)$ and $v_{n}\gg n^{-1/2}.$
Define $g_{n}$ and $\hat g_{n}$ as above with  ${\cal K}_{n}$ the set of all functions $g=(g(0), \ldots, g(S_{n}+1)) \in\RR^{S_{n}+2} $ such that $\Delta g(x)\geq 0$ for all $x\in\{1,\dots,S_{n}\}$ that satisfy $\Delta \hat p_n(x) \le v_n$. Then, $\hat g_{n}$ uniquely exists, both $\hat g_{n}$ and $\hat T_n:=\sum_{k=0}^{S_{n}+1}  \left(\hat g_n(k) - g_{n}(k)\right)^2$ are measurable, and conditionally on $X_{1},\dots,X_{n}$ we have  $\hat T_n
\stackrel{d}{\longrightarrow}\hat T_0$  in probability as $n\to\infty$, with $\hat T_0$ as in Theorem \ref{cor: asympt}.
\end{theorem}

\par \noindent We now state the main result of the section, again defining  $\hat T_{n}$  as in Theorem \ref{theo: CVcalib}.

\begin{theorem}\label{theo: test2}
Let $\alpha\in(0,1)$ and  $t_{\alpha,n}$  the conditional $(1-\alpha)$-quantile of  $\hat T_n$ given $X_{1},\dots,X_{n}$. If $p_0$ is convex on $\NN$ and supported on $\{0,\ldots, S\}$, then  $\limsup_{n \to \infty} P\big(T_n> {  t_{ \alpha,n}}\big) \le \alpha$.
\end{theorem}

In order to implement the test, we need to compute an approximation of 
 $t_{\alpha,n}$. This can be done using Monte-Carlo simulations as follows. Having observed $X_{1},\dots,X_{n}$, draw independent sequences $(Z_{i}^{(b)})_{1\leq i\leq S_{n}}$ for 
$b\in\{1,\dots,B\}$, where all variables $Z_{i}^{(b)}$ are i.i.d. standard Gaussian and $B>0$ is an integer. Then, for all $b$, compute
$g_{n}^{(b)}=\Gamma_{n}^{1/2}(Z_{0}^{(b)},\dots,Z_{S_{n}+1}^{(b)})^T$
and the least-squares projection $\hat g_{n}^{(b)}$ onto ${\cal K}_{n}$ using the algorithm described in \cite{BDK2014limit}. 
Then, $t_{\alpha,n}$ can be approximated by the $(1-\alpha)$-quantile of the empirical distribution
coresponding to $\sum_{k=0}^{S_{n}+1} (g_{n}^{(b)}(k)-\hat g_{n}^{(b)}(k))^2$, with $b\in\{1,\dots,B\}$.

\section{Calibrating under the least favorable hypothesis}
We consider below an alternative calibration that is easier to implement than the first one. Consider $\widetilde{\mathcal K}_0$ the set of convex functions $g$ on $\{0, \ldots, S+1 \}$; that is $g \in \tilde{\mathcal K}_0$ if and only if $\Delta g(k)\geq 0$ for all $k\in\{1,\dots,S\}$. 
Similarly, let $\widetilde{\mathcal K}_n$ be the set of convex functions on $\{0, \ldots, S_n +1 \}$. Let $\tilde{g}_0$ be the least squares projection of $g_0$ onto $\widetilde{\mathcal K}_0$ and $\tilde{g}_n$ that of $g_n$ onto $\widetilde{\mathcal K}_n$, with $g_0$ and $g_n$ as in Section \ref{sec: lim}. Finally, let $\tilde{t}_{n, \alpha}$ be the conditional $(1- \alpha)$-quantile of  $\widetilde T_n:=\sum_{k=0}^{S_{n}+1} ( \tilde{g}_n (k)- g_n (k))^2$ given  $(X_1, \ldots, X_n)$. Then, we have the following theorem. 

\begin{theorem}\label{theo: test3}
If  $p_0$ is convex on $\NN$ and supported on $\{0, \ldots, S\}$, then
$
\limsup_{n \to \infty} P\big(T_n> \tilde{t}_{n, \alpha}\big) \le \alpha
$ 
with equality if $p_0$ is the triangular pmf with support $\{0, \ldots, S \}$.
\end{theorem} 

The  test 
is asymptotically calibrated since the Type I error does not exceed $\alpha$. 
It reaches precisely $\alpha$ when $p_0$ is triangular, which 
can be viewed as the least favorable case for testing convexity. The theorem above does not exclude existence of other least favorable cases. 

\section{Simulations}

To illustrate the theory, we have considered four pmf's supported on $\{0, \ldots, 5 \}$. It follows from Theorem 7 in  \cite{durot2013least}  that any convex pmf 
on $\NN$ 
can be written as $\sum_{k \ge 1} \pi_k \mathcal{T}_j$ where $\pi_k \in [0,1], \sum_{k \ge 1}\pi_k =1$  and $\mathcal{T}_k(i) = 2 (k-i)_+ [k(k+1)]^{-1}$, the triangular pmf supported on $\{0, \ldots, k-1\}$. 
Under $H_0$, we considered the triangular pmf $p^{(1)}_0  =  \mathcal{T}_6$ 
and $p^{(2)}_0 = \sum_{k=1}^6 \pi_k \mathcal{T}_k$  with $\pi_1=0, \pi_2=\pi_3=1/6$, $\pi_4=0$ and $\pi_5 = \pi_6 =1/3$, which has knots at $2, 3$ and $5$.  Under $H_1$
we considered $p^{(1)}_1$ the pmf of a truncated Poisson on $\{0, \ldots, 5 \}$ with rate $\lambda =1.5$ and $p^{(2)}_1$ the pmf equal to $p^{(1)}_0$ on $\{2,\ldots, 5 \}$ such that $(p^{(2)}_1(0),p^{(2)}_1(1))=(p^{(1)}_0(0)  +0.008, p^{(1)}_0(1) -0.008)$. To investigate the asymptotic type I error and power of our tests, we have drawn 
$n \in \{500, 5000, 50000 \}$ rv's from the aforementioned pmf's. Here, $\alpha =5\%$ and
$\tilde{t}_{n, \alpha}$ 
was estimated for each drawn sample 
using $B=1000$ i.i.d.  copies of $g_n$. The rejection probability was estimated using $N=500$ replications of the whole procedure. For the first convexity test,  we considered the sequences $v_n \equiv \sqrt{\log (\log n)} n^{-1/2}$ and $n^{-1/4} $. We also added the sequence $v_n \equiv0$ to compare our approach with the naive one where no knot extraction is attempted.  The results are reported in Tables  \ref{PowerTestvn} and \ref{PowerLFH}.

\begin{table}
 \begin{tabular}{|c|ccc|ccc|ccc|}
    \hline
    \multirow{2}{*}{PMF} &
      \multicolumn{3}{c}{$n=500$} &
      \multicolumn{3}{c}{$n=5000$} &
      \multicolumn{3}{c|}{$n=50000$} \\
    & $0$ & $\frac{(\log(\log n))^{1/2}}{ n^{1/2}}$ & $n^{-1/4} $ & $0$ & $\frac{(\log(\log n))^{1/2}}{ n^{1/2}}$ &  $n^{-1/4}$&  $0$ & $\frac{(\log(\log n))^{1/2}}{ n^{1/2}}$  & $n^{-1/4}$ \\
    \hline
\hline
    $p^{(1)}_0$ & 0.226 & 0.106 & 0.054 & 0.286 & 0.092 & 0.062 & 0.256 & 0.086  & 0.050\\
    \hline
    $p^{(2)}_0$& 0.190 & 0.046 & 0.020 & 0.310 & 0.066& 0.018& 0.344 & 0.054  & 0.016\\
    \hline
    $p^{(1)}_1$ & 1 & 1 & 1 &1 & 1 & 1 &  1 & 1  & 1\\
    \hline
$p^{(2)}_1$ & 0.234 & 0.102 & 0.038 & 0.354 & 0.166 & 0.082 & 0.932 & 0.816  & 0.630\\
    \hline
\end{tabular}
\caption{Values of the asymptotic type I error for the pmfs $p^{(1)}_0$  and $p^{(2)}_0$ and power for $p^{(1)}_1$  and $p^{(2)}_1$  of the convexity test based on the sequence $(v_n)_n$.  The asymptotic level is $5\%$.}
  \label{PowerTestvn}
\end{table}

 \begin{table}
 \begin{tabular}{|c|c|c|c|}
    \hline
    PMF &  $n=500$    & $n=5000$  & $n=50000$ \\
      \hline
    $p^{(1)}_0$ & 0.048 & 0.044 & 0.058  \\
    \hline
    $p^{(2)}_0$&  0.014 & 0.032 & 0.020  \\
    \hline
    $p^{(1)}_1$ &  1 & 1 & 1 \\
    \hline
$p^{(2)}_1$ &  0.042 & 0.060 & 0.678  \\
    \hline
\end{tabular}
\caption{Values of the asymptotic type I error for the pmfs $p^{(1)}_0$  and $p^{(2)}_0$ and power for $p^{(1)}_1$  and $p^{(2)}_1$  of the test based on the least favorable hypothesis. The asymptotic level is $5\%$.}
  \label{PowerLFH}
\end{table}

The first conclusion is that  choosing $v_n =0$  does not give a valid test as the type I error can be as large as four times the targeted level! This can be explained by the fact that choosing $v_n=0$ makes the set on which $g_n$ is projected to be the largest possible and hence  $\Vert \hat{g}_n - g_n \Vert$ the smallest possible. This distance is hence stochastically smaller than the actual  limit, 
yielding a large probability of rejection.
The second conclusion is that the first test depends on the choice of $v_n$.
Small sequences  makes again the type I error large when the true convex pmf has only a tiny change in the slopes at its knots or has no knots as it is the case for $p^{(1)}_0  = T_6$. The question is then open as to how to choose such a sequence so that the test has the correct asymptotic level.  The second testing approach is, as expected, conservative when the true pmf is not triangular. For the true pmf  $p^{(1)}_1$ which strongly violates the convexity constraint, 
the power is equal to 1. 
For $p^{(2)}_1$,  which has only a small flaw at 2 with a change of slope equal to $-0.008$, the power values obtained with this second approach are comparable 
to those obtained with the first testing method and $v_n \equiv n^{-1/4}$. This is somehow expected as it is the largest sequence among the ones considered, yielding the largest distance between $g_n$ and its $L_2$ projection.

\section{Proofs}
In the sequel, for all $s>0$ and $u=(u(0), \ldots, u(s+1)) \in\RR^{s+2} $, we set $\|u\|_s=\sum_{k=0}^{s+1}(u(k))^2.$
\subsection{Preparatory lemmas}

\begin{lemma}\label{lem: mesurable}
Let $s>0$ be an integer, ${\cal K}\subset\RR^{s+2}$ a non-empty closed convex set, and $u=(u(0), \ldots, u(s+1)) \in\RR^{s+2} $. Then,
the minimizer of $\|g-u\|_{s}$ over $g\in{\cal K}$ uniquely exists. Moreover, denoting $\Phi(u)$ this minimizer,  the application $\Phi$ is measurable from $\RR^{s+2}$ to $\RR^{s+2}$; and the application $u \mapsto \|u-\Phi(u)\|_{s}$  is measurable. 
\end{lemma}

\par \noindent {\bf Proof:} It follows from standard results on convex optimization that $\Phi(u)$ uniquely exists for all $u$, and 
$\|\Phi(u)-\Phi(v)\|_{s}\leq \|u-v\|_{s}$ for all $u$ and $v$ in $\RR^{s+2}$.
This means  that $\Phi$ is a continuous function, whence it is measurable. Now,  the function  $u\mapsto(u,\Phi(u))$ is continuous, whence measurable. 
By continuity of the norm, this ensures that the application that maps $u$ into $\|u-\Phi(u)\|_{s}$ is measurable. \hfill{$\Box$}

\begin{lemma}\label{lem: inversible}
With $g_0$ as in Section \ref{sec: pb}, $(\Delta g_{0}(1),\dots,\Delta g_{0}(S))$  is a centered Gaussian vector with invertible dispersion matrix.
\end{lemma}

\par \noindent \textbf{Proof:}  For notational convenience, we assume in the sequel that $S\geq 3$. The case $S\leq 2$ can be handled likewise. Let $B$ be the $S\times(S+1)$-matrix which $j$-th line has components $j$,$j+1$ and $j+2$ equal to $1,-2$ and $1$ respectively while the other components are zero, for $j=1,\dots,S-1$, and $S$-th line has components equal to zero except the penultimate and the last one, which are equal respectively to $1$ and $-2$.
We have $p_{n}(S+1)=p_{0}(S+1)=0$ almost surely so that in the limit, $g_{0}(S+1)=0$ almost surely and
\begin{equation}\label{eq: vecDelta}
\left(
\Delta g_{0}(1),
\dots,
\Delta g_{0}(S)
\right)^T=B\left(
g_{0}(0),
\dots,
g_{0}(S)
\right)^T.\end{equation}
Hence, $(\Delta g_{0}(1),\dots,\Delta g_{0}(S))$  is a centered Gaussian vector with dispersion matrix 
\begin{equation}\label{eq: V}
V=B\Sigma_{0}B^T,
\end{equation}
where  $\Sigma_{0}$ is the dispersion matrix of the  vector on the right hand side of \eqref{eq: vecDelta}, {\it i.e.}  with component $(i+1,j+1)$ equal to $\mathbbm{1}_{\{i= j\}} p_{0}(i)-p_{0}(i)p_{0}(j)$  for all  $i,j=0,\dots,S$. Note that $\Sigma_{0}$ is obtained by deleting a line and a column of zeros in  $\Gamma_0$, the dispersion matrix of $g_0$.

It remains to prove that  $V$  is invertible. Let  $\sqrt p_{0}$ be the column vector in $\RR^{S+1}$ with components $\sqrt {p_{0}(0)}$, \dots, $\sqrt{p_{0}(S)}$ and let $\mbox{diag}(\sqrt p_{0})$ be the $(S+1)\times(S+1)$ diagonal matrix with diagonal components $\sqrt {p_{0}(0)},\dots,\sqrt{p_{0}(S)}$. 
Denoting by $I$ the identity matrix on $\RR^{S+1}$, 
the matrix (in the canonical basis) associated with the orthogonal projection from $\RR^{S+1}$ onto the orthogonal supplement of the linear space generated by $\sqrt p_{0}$ is given by
$I-\Pi_{0}=I-\sqrt p_{0} \sqrt p_{0}^T.$
The linear subspace of $\RR^{S+1}$ generated by $\sqrt p_{0}$ has dimension 1, so its  orthogonal supplement  in $\RR^{S+1}$ has dimension $S$, whence
$\mbox{rank}(I-\Pi_{0})=S.
$ 
Now,
$\Sigma_{0}=\mbox{diag}(\sqrt p_{0})(I-\Pi_{0})\mbox{diag}(\sqrt p_{0})$ where  diag$(\sqrt p_{0})$ is invertible, so we obtain 
\begin{equation}\label{eq: rank12}
\mbox{rank}(\Sigma_{0}^{1/2})=\mbox{rank}(\Sigma_{0})=S.
\end{equation}
This means that the kernel of $\Sigma_{0}^{1/2}$ is a linear subspace of $\RR^{S+1}$ whose dimension is equal to 1. Let us describe more precisely the kernel. Let $\lambda$ be the column vector in $\RR^{S+1}$ whose components are all equal to 1. Using that $\sum_{k=0}^Sp_{0}(k)=1$, it is easy to see that $\Sigma_{0}\lambda$ is the null vector in $\RR^{S+1}$. Therefore, $\lambda^T\Sigma_{0}\lambda=0$. This means that
$\|\Sigma_{0}^{1/2}\lambda\|^2
=0$
with $\|\;.\;\|$ the euclidian norm in $\RR^{S+1}$. Hence, $\Sigma_{0}^{1/2}\lambda$ is the null vector in $\RR^{S+1}$. This means that 
the kernel of $\Sigma_{0}^{1/2}$ is the linear subspace of $\RR^{S+1}$ generated by $\lambda$.

Next, let us determine the kernel of $V$. Let $\mu\in\RR^S$ with $V\mu=0$. Then, $\mu^TV\mu=0$ which, according to  \eqref{eq: V}, implies that
$\|\Sigma_{0}^{1/2}B^T\mu\|^2=\mu^TB\Sigma_{0}B^T\mu=0.$
This means that $\Sigma_{0}^{1/2}B^T\mu=0$. Since the kernel of $\Sigma_{0}^{1/2}$ is the linear subspace of $\RR^{S+1}$ generated by $\lambda$, we conclude that
$B^T\mu=a\lambda$
for some $a\in\RR$. Denote by $\mu_1,\dots,\mu_{S}$ the components of $\mu$. By definition of $B$ and $\lambda$, the vector $\mu$ satisfies the equations $
\mu_{1}=a$, $
\mu_{2}-2\mu_{1}=a$, $
\mu_{k-2}-2\mu_{k-1}+\mu_{k}=a \mbox{ for all }k\in\{3,\dots,S\}$ and $
\mu_{S-1}-2\mu_S=a.$
Arguing by induction, we obtain that this is equivalent to $
2\mu_{k}=ak(k+1) \mbox{ for all }k\in\{1,\dots,S\}$ and $
2\mu_S=\mu_{S-1}-a.$
Combining the first equation with $k=S,S-1$ to the second equation yields
$$aS(S+1)=-a+\frac{a(S-1)S}2.$$
Therefore,
$$a\left(1-\frac{(S-1)S}2+S(S+1)\right)=0.$$
This reduces to
 $a(2+S^2+3S)/2=0,$ which implies  that $a=0$. This mean that $\mu$ is the null vector in $\RR^{S}$ and therefore,  rank$(V)=S$. This means that $V$ is invertible.\hfill{$\Box$}

\begin{lemma}\label{lem: continu}
Let $\mathcal K _0$, $g_0$ and $\hat g_0$ be defined as in Section \ref{sec: pb}.  Then,  $\|\hat g_0-g_{0}\|_{S}$ has a continuous distribution.
\end{lemma}

\par \noindent \textbf{Proof:} \  Let  $F$ be the cumulative distribution function of $\|\hat g_{0}-g_{0}\|_{S}$. Note that this quantity is a properly defined random variable by the measurability proved in Lemma \ref{lem: mesurable}. 
We aim to prove that $F$ is a continuous function on $[0,\infty)$, using similar arguments as in the proof of Lemma 1.2 in \cite{gaenssler2007continuity}. 
First, we will prove that
\begin{equation}\label{eq: F>0}
F(t)>0\quad\mbox{for all }t\geq 0.
\end{equation}
To this end, note that 
$F(0)
=\PP(g_{0}\in \mathcal K_0). $
This means that
\begin{eqnarray}\label{eq: F(0)}
F(0)\geq\PP(\delta\in A)
\end{eqnarray}
with $\delta=(\Delta g_{0}(1),\dots,\Delta g_{0}(S))$ and $A$ the set of all vectors $(u_1,\dots,u_S)\in \RR^S$ such that $u_k\geq 0$ for all $k\in\{1,\dots,S\}$ with possible exceptions at points $k$ that are knots of $p_{0}$.  From Lemma \ref{lem: inversible}, the vector $\delta$ is a centered Gaussian vector whose dispersion matrix is invertible. Therefore, the vector possesses a density with respect to the Lebesgue measure on $\RR^S$ that is strictly positive on the whole space $\RR^S$. This implies that the probability that the vector belongs to a Borel set whose Lebesgue measure is not equal to zero, is strictly positive. In particular,  $\PP(\delta\in A)>0$. Thus, it follows from \eqref{eq: F(0)} that $F(0)>0.$
Combining this with the monotonicity of the function $F$ completes the proof of \eqref{eq: F>0}.

Next, we prove that the function $\log(F)$ (which is well defined on $[0,\infty)$ thanks to \eqref{eq: F>0}) is concave on $[0,\infty)$. For $u=(u(0), \ldots, u(S+1)) \in\RR^{S+2} $, let us write $\hat u$ the minimizer of 
$\sum_{k=0}^{S+1}  \left(g(k) - u(k)\right)^2$
over $g\in{\cal K}_{0}$.  For all $t\in[0,\infty)$, we define 
$A_{t}$ to be the set of all $u=(u(0),\dots,u(S+1))\in\RR^{S+2}$ such that $\|\hat u-u\|_{S}\leq t.$
Note that $A_{t}$ is a Borel set in ${\cal B}(\RR^{S+2})$ for all $t$ since, according to Lemma \ref{lem: mesurable}, the application $u\mapsto \|\hat u-u\|_{S}$ is  measurable.
Finally, let $\mu=\PP\circ g_{0}^{-1}$ be the distribution of $g_{0}$ on $\RR^{S+2}$ endowed with the Borel $\sigma$-algebra ${\cal B}(\RR^{S+2})$. This means that 
\begin{equation}\label{eq: mutoF}
F(t)=\mu(A_{t}).
\end{equation}

Fix $\lambda\in(0,1)$,  $t,t'\in[0,\infty)$, and consider an arbitrary $x\in\lambda A_{t}+(1-\lambda)A_{t'}$. Then, $x$ takes the form $x=\lambda u+(1-\lambda)v$ for some (non necessarily unique) $u\in A_{t}$ and $v\in A_{t'}$. By definition, both $\hat u$ and $\hat v$ belong to the convex set $\mathcal K_0$ and therefore, $\lambda \hat u+(1-\lambda)\hat v\in \mathcal K_0.$ Since $\hat x$ minimizes $\|g-x\|_{S}$
over  $g\in \mathcal K_0$, we conclude that
\begin{eqnarray*}
\|\hat x-x\|_{S}\leq\|\lambda \hat u+(1-\lambda)\hat v-x\|_S\leq \|\lambda (\hat u-u)+(1-\lambda)(\hat v-v)\|_S,
\end{eqnarray*}
using that $x=\lambda u+(1-\lambda)v$. It then follows from the triangle inequality that 
\begin{eqnarray*}
\|\hat x-x\|_{S}
&\leq&\lambda\|\hat u-u\|_S+(1-\lambda)\|\hat v-v\|_S.
\end{eqnarray*}
Since $u\in A_{t}$ and $v\in A_{t'}$, we have $\|\hat u-u\|_{S}\leq t $ and $\|\hat v-v\|_{S}\leq t' $, which implies that 
$
\|\hat x-x\|_{S}\leq\lambda t+(1-\lambda)t'.
$ 
Hence, $x\in A_{\lambda t+(1-\lambda)t'}.$ This means that
\begin{eqnarray}\label{eq: inclusionA}
\lambda A_{t}+(1-\lambda)A_{t'}\subset A_{\lambda t+(1-\lambda)t'}.
\end{eqnarray}
Now, $\mu=\PP\circ g_{0}^{-1}$ is a Gaussian probability measure on ${\cal B}(\RR^{S+2})$, so it follows from Lemma 1.1 in \cite{gaenssler2007continuity} that $\mu$ is log-concave in the sense that
$$\mu_{\star}(\lambda A+(1-\lambda)B)\geq \mu(A)^{\lambda}\mu(B)^{1-\lambda}$$
for all $\lambda\in(0,1)$ and $A,B\in {\cal B}(\RR^{S+2})$, with $\mu_{\star}$  the inner measure pertaining to $\mu$.  Applying this with $A=A_t$ and $B=A_{t'}$, and combining  
with \eqref{eq: inclusionA} yields
$$\mu_{\star}(A_{\lambda t+(1-\lambda)t'})\geq \mu(A_{t})^{\lambda}\mu(A_{t'})^{1-\lambda}$$
for all $t,t'\in[0,\infty)$.
The same inequality remains true with $\mu_{\star}$ replaced by $\mu$.
 Using \eqref{eq: mutoF}, and taking the logarithm on both sides of the inequality, we conclude that
$$\log\big(F(\lambda t+(1-\lambda)t')\big)\geq \lambda\log\big(F(t)\big)+(1-\lambda)\log\big(F(t')\big)$$
for all $\lambda\in(0,1)$, and $t,t'\in[0,\infty)$. This means that the function $\log(F)$ is concave on $[0,\infty)$. Recalling \eqref{eq: F>0}, we conclude that the function $\log(F)$ is continuous on $[0,\infty)$, whence  $F$ is continuous on $[0,\infty)$. This completes the proof of Lemma \ref{lem: continu}. \hfill$\Box$

\subsection{Proofs of the main results}

\par \noindent \textbf{Proof of Theorem \ref{cor: asympt}:}
Assume that $p_0$ is  convex on $\NN$ and supported on $\{0,\ldots, S\}$. It can be proved  in the same manner as Theorem 3.2 in \cite{BDK2014limit} that 
\begin{eqnarray}\label{Asymp}
\sqrt n \left (\widehat p_n - p_0, p_{n}-p_{0}\right) \Rightarrow
\left ( 
\widehat g_{0},
g_{0}
\right)\mbox{ as }n \to \infty. 
\end{eqnarray}
as a joint weak convergence on $\{0,\dots,S+1\}$. Now, it follows from  \citet[Proposition 3.5]{BDK2014limit} that with probability one, $\hat p_{n}$ is supported on $\{0,\dots,S+1\}$ for sufficiently large $n$,  and $\hat p_{n}$ also is supported on that set by definition.
%
%
%
Hence, 
$\| p_{n}-\hat p_{n}\|=\| p_{n}-\hat p_{n}\|_{S}.$
Combining this with  \eqref{Asymp} completes the proof of the theorem. \hfill $\Box$


\par \noindent {\bf Proof of Theorem \ref{theo: CVcalib}:}  \  Clearly, ${\cal K}_{n}$ is a non-empty closed convex subset of $\RR^{S_{n}+2}$.
Hence, we consider the specific case of $s=S_{n}$ and ${\cal K}={\cal K}_{n}$ in Lemma \ref{lem: mesurable}. In the notation of the lemma, we have $\hat g_{n}:=\Phi(g_{n})$, so that $\hat g_{n}$ is uniquely defined. Moreover, since both $\Phi$ and $g_{n}$ are measurable, we conclude that $\hat g_n=\Phi(g_{n})$ is measurable.  Likewise, $\|g_{n}-\hat g_{n}\|_{S_{n}}$  is measurable. This proves the fist two assertions in Theorem \ref{theo: CVcalib}. Next, similar to \citet[Theorem 3.3]{BDK2014limit}, the following joint weak convergence on $\{0,\dots,S+1\}$ can be proved: conditionally on $X_{1},\dots,X_{n}$,
$
\left (
\widehat g_n  ,
g_{n}
\right) \Rightarrow
\left ( 
\widehat g_{0},
g_{0}
\right)\mbox{ in probability as }n \to \infty. 
$ 
The result follows, since $S_n=S$ with provability that tends to one.
%
%
%
 \hfill{$\Box$}


\par \noindent \textbf{Proof of Theorem \ref{theo: test2}:}  \ Assume that $p_{0}$ is convex on $\NN$ with support $\{0,\dots,S\}$.
By Theorem \ref{cor: asympt},   $T_{n}$ converges in distribution to $\|\hat g_{0} - g_{0}\|_{S}$ as $n\to\infty$. Combining this with Lemma \ref{lem: continu} together with the fact that convergence in distribution to a continuous distribution implies uniform convergence of the corresponding distribution functions yields
\begin{eqnarray*}
P\left(T_n>t_{\alpha,n}\right)&=&P\left(\|\hat g_{0}-g_{0}\|_{S}>t_{\alpha,n}\right)+o(1).
\end{eqnarray*}
Now, it follows from Theorem \ref{theo: CVcalib} that
\begin{eqnarray*}
P\left(\|\hat g_{0}-g_{0}\|_{S}>t_{\alpha,n}\right)&=&P\left(\|\hat g_{n}-g_{n}\|_{S_{n}}>t_{\alpha,n}\ |\ X_{1},\dots,X_{n}\right)+o_{p}(1)
\end{eqnarray*}
where by definition of $t_{\alpha,n}$, the probability on the right-hand side is less than or equal to $\alpha$
for all $n$. Combining this with the preceding  display completes the proof. \hfill{$\Box$}


\par \noindent \textbf{Proof of Theorem \ref{theo: test3}:} \ We have $ \|\tilde{g}_{0}-g_{0}\|_{S} \ge \|\hat g_{0}-g_{0}\|_{S}$ since $\widetilde{\mathcal K}_n\subset {\mathcal K}_0$, whence
\begin{eqnarray}\label{MajProb}
P\left(T_n> \tilde{t}_{\alpha,n}\right)&=& P\left(\|\hat g_{0}-g_{0}\|_{S}> \tilde{t}_{\alpha,n}\right)+o(1) \nonumber\\
                                                                                 & \le & P\left(\|\tilde{g}_{0}-g_{0}\|_{S}> \tilde{t}_{\alpha,n}\right)+o(1),
\end{eqnarray}
using similar arguments as for the proof of Theorem \ref{theo: test2} for the equality.
Using arguments similar to those in Theorem 3.3 of \cite{BDK2014limit}, we can  show that conditionally on $X_{1},\dots,X_{n}$,
$
\left (
\tilde{g}_n    ,
g_{n}
\right)\to
\left ( 
\tilde{g}_{0},
g_{0}
\right)\mbox{almost surely as }n \to \infty. 
$ 
It follows that conditionally on $X_{1},\dots,X_{n}$, $\Vert \tilde{g}_n - g_n \Vert_{S_n} $ converges weakly to $\Vert \tilde{g}_0 - g_0 \Vert_S$ almost surely whence
\begin{eqnarray*}
P\left(T_n> \tilde{t}_{\alpha,n}\right) & \le &  P\left(\|\tilde{g}_n-g_n\|_{S_n}> \tilde{t}_{\alpha,n} | X_1, \ldots, X_n \right)+o(1) \\
& = & \alpha + o(1), \ \textrm{by definition of $\tilde{t}_{\alpha, n}$.}
\end{eqnarray*}
The inequality in (\ref{MajProb}) becomes an equality if $p_0$ is triangular, so the theorem follows. \hfill $\Box$

\bibliographystyle{ims}
\bibliography{TestConvex}

\end{document}